\documentclass[12pt]{amsart}
\usepackage{epsfig}
\begin{document}
\title{Toric varieties whose blow-up at a point is Fano}
\author{Laurent BONAVERO}
\maketitle
\noindent
\def\restriction{\string |}
\newcommand{\pp}{\rm ppcm}
\newcommand{\pg}{\rm pgcd}
\newcommand{\Ker}{\rm Ker}
\newcommand{\C}{{\mathbb C}}
\newcommand{\Q}{{\mathbb Q}}
\newcommand{\B}{\rm B}
\newcommand{\GL}{\rm GL}
\newcommand{\SL}{\rm SL}
\newcommand{\diag}{\rm diag}

\def\refname{References}

\newcommand{\finpreuve}{{\nobreak\hfil\penalty50\hskip2em
\hbox{}\nobreak\hfil\mbox{\rule{1ex}
{1ex}\qquad}\parfillskip=-20pt\finalhyphendemerits=0\par\medskip}}

\newtheorem{theo}{Theorem}
\newtheorem{prop}{Proposition}
\newtheorem{lemm}{Lemma}
\newtheorem{lemmf}{Lemme fondamental}
\newtheorem{defi}{D\'efinition}
\newtheorem{exo}{Exercice}
\newtheorem{rem}{Remarque}
\newtheorem{cor}{Corollaire}
\newcommand{\CC}{{\mathbb C}}
\newcommand{\ZZ}{{\mathbb Z}}
\newcommand{\RR}{{\mathbb R}}
\newcommand{\QQ}{{\mathbb Q}}
\newcommand{\FF}{{\mathbb F}}
\newcommand{\PP}{{\mathbb P}}
\newcommand{\codim}{\operatorname{codim}}
\newcommand{\Ho}{\operatorname{Hom}}
\newcommand{\Pic}{\operatorname{Pic}}
\newcommand{\NE}{\operatorname{NE}}
\newcommand{\Nun}{\operatorname{N}}
\newcommand{\card}{\operatorname{card}}
\newcommand{\Hilb}{\operatorname{Hilb}}
\newcommand{\mult}{\operatorname{mult}}
\newcommand{\vol}{\operatorname{vol}}
\newcommand{\divi}{\operatorname{div}}
\newcommand{\pr}{\operatorname{pr}}
\newcommand{\con}{\operatorname{cont}}
\newcommand{\ima}{\operatorname{Im}}

\newcounter{subsub}[subsection]
\def\thesubsub{\thesubsection .\arabic{subsub}}
\def\subsub#1{\addtocounter{subsub}{1}\par\vspace{3mm}
\noindent{\bf \thesubsub ~ #1 }\par\vspace{2mm}}
\def\coker{\mathop{\rm coker}\nolimits}
\def\pr{\mathop{\rm pr}\nolimits}
\def\im{\mathop{\rm Im}\nolimits}
\def\hfl#1#2{\smash{\mathop{\hbox to 12mm{\rightarrowfill}}
\limits^{\scriptstyle#1}_{\scriptstyle#2}}}
\def\vfl#1#2{\llap{$\scriptstyle #1$}\big\downarrow
\big\uparrow
\rlap{$\scriptstyle #2$}}
\def\diagram#1{\def\normalbaselines{\baselineskip=0pt
\lineskip=10pt\lineskiplimit=1pt}   \matrix{#1}}
\def\limind{\mathop{\oalign{lim\cr
\hidewidth$\longrightarrow$\hidewidth\cr}}}

\long\def\InsertFig#1 #2 #3 #4\EndFig{
\hbox{\hskip #1 mm$\vbox to #2 mm{\vfil\includegraphics{#3}}#4$}}
\long\def\LabelTeX#1 #2 #3\ELTX{\rlap{\kern#1mm\raise#2mm\hbox{#3}}}

{\let\thefootnote\relax
}

{\bf Abstract.} We classify smooth toric Fano varieties of
dimension $n\geq 3$ containing a toric divisor isomorphic to the
$(n-1)$-dimensional projective space. As a consequence of this
classification, we show that any smooth complete toric variety
$X$ of dimension $n\geq 3$ with a fixed point $x\in X$ such that
the blow-up $B_x(X)$ of $X$ at $x$ is Fano is isomorphic either
to the $n$-dimensional projective space or to the blow-up of the
$n$-dimensional projective space along an invariant linear
codimension two subspace. As expected, such results are proved
using toric Mori theory due to Reid.

\section*{Introduction}

\hspace{3mm} Smooth blow-ups and blow-downs between toric smooth
Fano varieties have been intensively studied: see \cite{Bat82},
\cite{Bat99}, \cite{Oda88}, \cite{WWa82} and more recently
\cite{Sat00}. In this Note, we prove the following result using
toric Mori theory (see also the uncorrect exercise V.3.7.10
mentionned in \cite{Kol99}): As usual, $T$ denotes the big torus
acting on a toric variety; if $Y$ is a smooth subvariety of a
smooth variety $X$, $B_Y(X)$ denotes the blow-up of $X$ along $Y$
and a variety $X$ is called Fano if and only if $-K_X$ is ample.

\begin{theo}\label{fano}
Let $X$ be a smooth and complete toric variety of dimension $n\geq 3$.
Suppose there exists a $T$-fixed point $x$ in $X$ such
that $B_x(X)$ is Fano.
Then either $X \simeq \PP^n$ and $x$ can be choosen arbitrary
or $X \simeq B_{\PP^{n-2}}(\PP ^n)$ and $x$
must be choosen outside the exceptional divisor.
\end{theo}

\hspace{3mm} Let us say that when $X$ is a toric surface with a
$T$-fixed point $x$ in $X$ such that $B_x(X)$ is Fano, then $X$ is
isomorphic to $\PP ^2$ blown-up at $m$ $T$-fixed points with
$m=0,1$ or $2$ or to $\PP ^1 \times \PP ^1$  blown-up at $m$
$T$-fixed points with $m=0$ or $1$. Recall also that smooth Fano
toric varieties are classified in dimension less or equal to $4$
(\cite{Bat82}, \cite{Bat99}, \cite{Oda88}, \cite{WWa82} and
\cite{Sat00}) together with smooth blow-ups and blow-downs between
them; in particular, Theorem \ref{fano} could be proved in
dimension $3$ and $4$ just by looking at the classification.

\hspace{3mm} In fact, we will obtain Theorem \ref{fano} as a
consequence of the following result (which is inspired by a
private communication of J.~Wi\'sniewski):

\begin{theo}\label{divisor}
Let $X$ be a smooth toric Fano variety of dimension $n\geq 3$.
Then, there exists a toric divisor $D$ of $X$ isomorphic to
$\PP^{n-1}$ with ${\mathcal O}_{\PP ^{n-1}}(d)$ as normal bundle
in $X$ if and only if one of the following situations occurs{\rm :}
\begin{enumerate}
\item [(i)]
$X \simeq \PP^n$, $d=1$ and $D$ is a linear codimension one subspace
of $X$,
\item [(ii)]
$X \simeq \PP ({\mathcal O}_{\PP ^{1}} \oplus
{\mathcal O}_{\PP ^{1}}(1)^{\oplus n-1}) \simeq B_{\PP^{n-2}}(\PP^n)$,
$d=0$ and $D$ is a fiber of the projection on $\PP^1$,
\item [(iii)] there exists an integer $\nu$ satisfying $0\leq \nu\leq n-1$
such that $X$ is isomorphic to $\PP ({\mathcal O}_{\PP ^{n-1}}
\oplus {\mathcal O}_{\PP^{n-1}}(\nu))$ and $D$ is either the
divisor $\PP ({\mathcal O}_{\PP ^{n-1}})$ (and $d=\nu$) or the
divisor $\PP ({\mathcal O}_{\PP ^{n-1}}(\nu))$ (and $d=-\nu$),
\item [(iv)] there exists an integer $\nu$ satisfying $0\leq \nu \leq n-2$
such that $X$ is isomorphic to the blow-up of $\PP ({\mathcal
O}_{\PP ^{n-1}} \oplus {\mathcal O}_{\PP ^{n-1}}(\nu+1))$ along a
linear $\PP ^{n-2}$ contained in the divisor $\PP ({\mathcal
O}_{\PP ^{n-1}})$ and $D$ is either the strict transform of the
divisor $\PP ({\mathcal O}_{\PP ^{n-1}})$ (and $d=\nu$) or the
strict transform of the divisor $\PP ({\mathcal O}_{\PP
^{n-1}}(\nu +1))$ (and $d=-\nu-1$).
\end{enumerate}
\end{theo}

\hspace{3mm} Remark that the adjunction formula implies that $d
\geq 1-n$. As an immediate consequence of Theorem \ref{divisor},
{\em there are exactly $2n+1$ distinct smooth toric Fano varieties
of dimension $n\geq 3$ containing a toric divisor isomorphic to
$\PP^{n-1}$.}

\section{Notation}

\hspace{3mm} We briefly review notation and very basic facts of
toric geometry (see \cite{Ful93} or \cite{Oda88} for details).

\smallskip

\hspace{3mm} A toric variety $X$ is defined by a fan $\Delta$ in a
lattice $N$ (the elements of $N$ are the one parameter subgroups
of the big torus $T$). If $X$ is smooth, any cone of $\Delta$ is
simplicial, generated by a family of lattice vectors which is
part of a basis of $N$. Any such cone $\langle e_1,\ldots,e_r
\rangle$ defines a smooth $T$-invariant subvariety of codimension
$r$ which is the closure of a $T$-orbit.
Recall that on a toric variety $X$, a $T$-invariant
Cartier divisor is ample if and only if its
intersection with any toric curve of $X$ is strictly positive
\cite{Oda88}.

\smallskip

\hspace{3mm} The cone of effective curves modulo numerical
equivalence (usually denoted by $\NE (X)$) of a smooth projective
toric variety is polyhedral generated by the $T$-invariant curves
of $X$ \cite{Rei83}. A $T$-invariant extremal curve $C$ of $X$ is
called Mori extremal if moreover $-K_X \cdot C >0$. Finally, if
$C$ is a $T$-invariant extremal curve with normal bundle $N_{C/X}=
\bigoplus _{i=1}^{n-1}{\mathcal O}_{\PP ^1}(a_i)$ generating an
extremal ray $R$ of $\NE (X)$, let $$\alpha = \card \{ i \in
[1,\ldots,n-1] \, |\, a_i <0\} \,\,  \mbox{and} \, \, \beta =
\card \{ i \in [1,\ldots,n-1] \, |\, a_i \leq 0\}.$$ Then, toric
Mori theory, due to Reid \cite{Rei83}, says that the contraction
of $R$ defines a map $\varphi_R: X \to Y$, which is birational if
and only if $\alpha \neq 0$. In that case, its exceptional locus
$A(R)$ in $X$ is $(n-\alpha)$-dimensional, $B(R)=\varphi_R (A(R))$
is $(\beta-\alpha)$-dimensional and the restriction of
$\varphi_R$ to $A(R)$ is a flat morphism, with fibers isomorphic
to weighted projective spaces. If $\alpha =0$, $\varphi_R: X \to
Y$ is a smooth $\PP ^{n-\beta}$-fibration and $Y$ is smooth and
projective.

\section{Fano varieties with a divisor isomorphic to a projective space.}

\hspace{3mm} In this section, we prove Theorem \ref{divisor}.

\subsection{Mori contraction on $X$}

In this subsection,
$X$ is a smooth toric Fano variety of dimension $n\geq 3$
containing a toric divisor $D$ isomorphic
to $\PP^{n-1}$ and $N_{D/X}={\mathcal O}_{\PP ^{n-1}}(d)$.
Let $[l_D] \in \NE (X)$ be the class in $\NE (X)$
of a line $l_D$ contained in $D$ (this class does not depend
on the choice of the line).

\begin{prop}\label{transverse} Suppose
there exists a Mori extremal curve $\omega$ transverse to $D$
such that $[\omega] \in \NE (X)$ does not belong to the ray
generated by $[l_D]$. Denote by $\varphi _{[\omega]}$ the Mori
contraction defined by $\omega$. Then
\begin{enumerate}
\item [(i)] either $\nu := |d|$ satisfies $0 \leq \nu \leq n-1$,
$X \simeq \PP ({\mathcal O}_{\PP ^{n-1}} \oplus {\mathcal O}_{\PP
^{n-1}}(\nu))$ and $\varphi _{[\omega]}~: X \to \PP ^{n-1}$ is the
natural fibration, or
\item [(ii)] there exists a smooth toric Fano variety $X'$ with a
$T$-invariant smooth divisor $D'$ such that $$(D',N_{D'/X'})
\simeq (\PP ^{n-1},{\mathcal O}_{\PP ^{n-1}}(d+1)),$$ $\varphi
_{[\omega]}: X  \to X'$ is the blow-up of $X'$ along a toric
subvariety $Y \simeq \PP ^{n-2}$ contained in $D'$ and $D$ is the
strict transform of $D'$.
\end{enumerate}
\end{prop}

In Case (ii), we get a new smooth toric Fano variety $X'$
containing a toric divisor $D'$ isomorphic to $\PP ^{n-1}$. This
motivates the following definition.

\medskip

\noindent {\bf Definition.} When the situation (ii) of
Proposition {\rm \ref{transverse}} occurs, we say that the pair
$(X,D)$ {\em can be simplified}.

\medskip

\noindent {\it Proof of Proposition} {\rm \ref{transverse}.} Let
$$\displaystyle{N_{\omega/X} = \bigoplus _{i=1}^{n-1}{\mathcal
O}_{\PP ^1}(a_i)}$$ be the normal bundle of $\omega$ in $X$ and
as in the previous section: $$\alpha = \card \{ i \in
[1,\ldots,n-1] \, |\, a_i <0\} \,\,  \mbox{and} \, \, \beta =
\card \{ i \in [1,\ldots,n-1] \, |\, a_i \leq 0\}.$$ Since
$[\omega] \in \NE (X)$ does not belong to the ray generated by
$[l_D]$, each $a_i$ is less or equal to zero. Therefore, since
$\displaystyle{-K_X \cdot \omega = 2 + \sum_{i=1}^{n-1}a_i > 0}$,
there are only two possibilities:
\begin{enumerate}
\item [(i)] every $a_i =0$; therefore $\alpha =0$, $\beta = n-1$
and the Mori contraction $\varphi _{[\omega]}: X\to Z$ is a $\PP
^1$-fibration on $Z$. Since $D \simeq \PP ^{n-1}$ is a section of
this fibration (by the transversality assumption, $D\cdot \omega
=1$), we get $Z \simeq \PP ^{n-1}$ and, if $\nu := |d|$, $X$ is
isomorphic to $\displaystyle{\PP ({\mathcal O}_{\PP ^{n-1}} \oplus
{\mathcal O}_{\PP ^{n-1}}(\nu))}$ which is Fano if and only if $0
\leq \nu \leq n-1$, or
\item [(ii)] there is exacly one of the
$a_i$'s equal to $-1$ and each other equal to $0$. Therefore
$\alpha =1$, $\beta = n-1$ and $\varphi _{[\omega]}:X\to X'$ is a
smooth blow-down on a smooth codimension two center. Denote by $E
\subset X$ the exceptional divisor of $\varphi _{[\omega]}$.
Since $D \cdot \omega =1$, the center of the blow-up is
isomorphic to $E \cap D$, i.e., isomorphic to $\PP ^{n-2}$.
Therefore, since $N_{D/X}= {\mathcal O}_{\PP ^{n-1}}(d)$, the
center of the blow-up $\varphi _{[\omega]}$ in $X'$ is isomorphic
to $\PP ^{n-2}$ with normal bundle ${\mathcal O}_{\PP
^{n-2}}(d+1)\oplus {\mathcal O}_{\PP ^{n-2}}(1)$. Therefore $X'$
is Fano by Lemma \ref{blow} below (recall that $d\geq 1-n$).
Moreover, $D':=\varphi _{[\omega]} (D)$ is a $T$-invariant smooth
divisor containing the center of the blow-up $\varphi
_{[\omega]}$ and satisfying $\displaystyle{(D',N_{D'/X'}) \simeq
(\PP ^{n-1},{\mathcal O}_{\PP ^{n-1}}(d+1)).}$\finpreuve
\end{enumerate}

\begin{lemm}\label{blow}
Let $X$ be a smooth toric variety of dimension $n$. Suppose there
exists a $T$-invariant subvariety $Y$ isomorphic
to $\PP ^{n-2}$ with normal bundle
${\mathcal O}_{\PP ^{n-2}}(a) \oplus {\mathcal O}_{\PP ^{n-2}}(b)$
such that $B_Y(X)$ is Fano.
Then $X$ is Fano if and only if $n-1+a+b >0$.
\end{lemm}

\noindent {\it Proof.} Since $B_Y(X)$ is Fano, $-K_X$ has
strictly positive intersection with any curve not contained in
$Y$, and if $C$ is a line contained in $Y$, then $-K_X \cdot C =
n-1+a+b$. \finpreuve

\medskip

\hspace{3mm} Let us end this part by the following lemma, which
says that Case (ii) in Proposition \ref{transverse} can not occur
twice consecutively:

\begin{lemm}\label{twice}
With the previous notation, assume that the pair $(X,D)$ can be
simplified, and let $\varphi_{[\omega]}:X\to X'$ be the
corresponding codimension two smooth blow-down as in Proposition
{\rm \ref{transverse} (ii)}. Then, the pair
$(X',\varphi_{[\omega]}(D))$ can not be simplified.
\end{lemm}

\noindent {\it Proof.} By contradiction, suppose
$(X',\varphi_{[\omega]}(D))$ can be simplified and denote by
$\varphi_{[\omega ']}: X' \to X''$ the corresponding smooth
codimension two blow-down. The exceptional divisor $E' \subset X'$
of $\varphi_{[\omega ']}$ intersects $D':=\varphi_{[\omega]}(D)$
along a $\PP^{n-2}$ which itself meets the center $Z$ of
$\varphi_{[\omega]}$ (since two $\PP ^{n-2}$ contained in a
$\PP^{n-1}$ must intersect). Let $C$ be the fiber of
$\varphi_{[\omega ']}$ containing a given point of $E' \cap D'
\cap Z$. Then the strict transform of $C$ in $X$ is a curve with
normal bundle in $X$ equals to ${\mathcal O}_{\PP ^{1}}^{\oplus
n-2}\oplus {\mathcal O}_{\PP ^{1}}(-2)$, and hence with zero
intersection on $-K_X$, a contradiction, since $X$ is
Fano.\finpreuve

\subsection{Proof of Theorem \ref{divisor}}

As before, $X$ is a smooth toric Fano variety of dimension $n\geq 3$
containing a toric divisor $D$ isomorphic
to $\PP^{n-1}$ and $N_{D/X}={\mathcal O}_{\PP ^{n-1}}(d)$.
Let $[l_D] \in \NE (X)$ be the class in $\NE (X)$
of a line contained in $D$.

\begin{prop}\label{d positif}
Suppose that $d \geq 0$ and $[l_D]$ spans an extremal ray of $\NE
(X)$. Then
\begin{enumerate}
\item [(i)] either $d =0$ and $X \simeq \PP^1 \times \PP^{n-1}$
or $X \simeq \PP ({\mathcal O}_{\PP ^{1}} \oplus {\mathcal
O}_{\PP ^{1}}(1)^{\oplus n-1})$, or
\item [(ii)] $d=1$ and $X\simeq \PP ^n$.
\end{enumerate}
\end{prop}

\noindent {\it Proof.} If $l_D$ is a line contained in $D$, then
$$ N_{l_D/X} = {\mathcal O}_{\PP ^1}(1)^ {\oplus n-2} \oplus
{\mathcal O}_{\PP ^1}(d).$$ If $d=0$, then the Mori contraction
$\varphi_{[l_D]}$ is a smooth $\PP ^{n-1}$-fibration on $\PP ^1$,
therefore $X$ is isomorphic to $\PP (\bigoplus_{i=1}^{n}{\mathcal
O}_{\PP ^1}(a_i))$, which is Fano if and only if $X \simeq \PP^1
\times \PP^{n-1}$ or $X \simeq \PP ({\mathcal O}_{\PP ^{1}} \oplus
{\mathcal O}_{\PP ^{1}}(1)^{\oplus n-1})$. If $d > 0$, the Mori
contraction $\varphi_{[l_D]}$ maps $X$ to a point, therefore
$X\simeq \PP ^n$ and $d=1$.\finpreuve

\medskip

\hspace{3mm} Now, we are ready to prove Theorem \ref{divisor}: Let
$X$ be a smooth toric Fano variety of dimension $n\geq 3$. Suppose
there exists a toric divisor $D$ of $X$ isomorphic to $\PP^{n-1}$,
and let ${\mathcal O}_{\PP ^{n-1}}(d)$ be its normal bundle in
$X$. Let also $[l_D] \in \NE (X)$ be the class in $\NE (X)$ of a
line contained in $D$.

\begin{enumerate}
\item [$\bullet$] First case: Suppose that either $d <0$
or $d \geq 0$ and $[l_D]$ does not span an extremal ray in $\NE
(X)$. Since $D$ is effective, there exists a Mori extremal curve
$\omega$ transverse to $D$ such that $[\omega] \in \NE (X)$ does
not belong to the ray generated by $[l_D]$. Therefore Proposition
\ref{transverse} applies: $X \simeq \PP ({\mathcal O}_{\PP ^{n-1}}
\oplus {\mathcal O}_{\PP ^{n-1}}(|d|))$ (and $0 < |d| \leq n-1$)
or the pair $(X,D)$ can be simplified.
\item [$\bullet$] Second case: $d \geq 0$ and
$[l_D]$ spans an extremal ray of $\NE (X)$. Then apply
Proposition~\ref{d positif}.
\end{enumerate}
\hspace{3mm} As a result, either $X$ satisfies one of the
conclusions (i), (ii) or (iii) of Theorem \ref{divisor}, or the
pair $(X,D)$ can be simplified. In the latter case, let
$\varphi_{[\omega]}:X\to X'$ be the corresponding codimension two
smooth blow-down as in Proposition \ref{transverse} (ii). Then,
since the pair $(X',D')$ can not be simplified by Lemma
\ref{twice}, applying the same process to the Fano variety $X'$
with $D' = \varphi_{[\omega]}(D)$ and $d'=d+1$, $X'$ itself must
satisfy one of the conclusions (i), (ii) or (iii) of Theorem
\ref{divisor}. In case $X'$ is isomorphic to $\PP ^n$, we get
that $X \simeq \PP ({\mathcal O}_{\PP ^{1}} \oplus {\mathcal
O}_{\PP ^{1}}(1)^{\oplus n-1})$. Moreover, $X'$ can not be
isomorphic to $\PP ({\mathcal O}_{\PP ^{1}} \oplus {\mathcal
O}_{\PP ^{1}}(1)^{\oplus n-1})$, because assuming the contrary,
$(X,D)$ could be simplified twice, a contradiction with Lemma
\ref{twice}. Finally, suppose $X' \simeq \PP ({\mathcal O}_{\PP
^{n-1}} \oplus {\mathcal O}_{\PP ^{n-1}}(|d+1|))$. Since $X'$ is
Fano, we get $0 \leq |d+1| \leq n-1$, which together with the
inequality $d \geq 1-n$ shows that $X$ satisfies conclusion (iv)
of Theorem \ref{divisor}.\finpreuve

\section{Proof of Theorem~1}

\hspace{3mm} Let $X$ be a smooth toric complete variety of
dimension $n\geq 3$. Suppose in the sequel that there exists a
$T$-fixed point $x$ in $X$ such that $B_x(X)$ is Fano (it is
well-known that $X$ is therefore also Fano). Hence $B_x(X)$ is a
Fano variety containing a toric divisor (the exceptional divisor
of the blow-up $\pi : B_x(X) \to X$) isomorphic to $\PP ^{n-1}$
with normal bundle ${\mathcal O}_{\PP ^{n-1}}(-1)$. Applying
Theorem \ref{divisor} to $B_x(X)$ with $d=-1$ gives that either
\begin{enumerate}
\item[$\bullet$] $B_x(X) \simeq \PP ({\mathcal O}_{\PP ^{n-1}} \oplus {\mathcal
O}_{\PP ^{n-1}}(-1))$ therefore $X \simeq \PP^n$, or
\item[$\bullet$] $B_x(X)$ is isomorphic to
the blow-up of $\PP ({\mathcal O}_{\PP ^{n-1}} \oplus {\mathcal
O}_{\PP ^{n-1}}) = \PP^1 \times \PP^{n-1}$ along a $\PP ^{n-2}$
contained in a fiber of the projection $\PP^1 \times \PP^{n-1}
\to \PP^1$. Therefore, $X \simeq \PP ({\mathcal O}_{\PP ^{1}}
\oplus {\mathcal O}_{\PP ^{1}}(1)^{\oplus n-1}) \simeq
B_{\PP^{n-2}}(\PP^n)$ and $x$ is outside the exceptional divisor
of the blow-up $B_{\PP^{n-2}}(\PP^n)\to \PP^n$.\finpreuve
\end{enumerate}

-----------

{\em \noindent Institut Fourier, UMR 5582,
Universit\'e de Grenoble 1, BP 74. 38402 SAINT-MARTIN d'H\`ERES. FRANCE. \\
\noindent e-mail: bonavero@ujf-grenoble.fr }

\end{document}